\documentclass[11pt,a4paper,reqno]{amsart}
\usepackage{stmaryrd}
\usepackage[english]{babel}
\usepackage{pst-grad} 
\usepackage{pst-plot} 
\usepackage{pstricks}
\usepackage{amsmath,amssymb,mathrsfs,amsthm, mathtools}
\usepackage{tikz} 
\usepackage{bbm}
\usepackage{mathcomp,wasysym}  
\usepackage{graphicx}  
\usepackage[all]{xy} \xyoption{arc} \xyoption{color}
\usepackage{epsfig}
\usepackage{cite}
\newcommand{\pare}[1]{\left( #1 \right)}
\newcommand{\system}[1]{\left\lbrace #1 \right.}
\newcommand{\norm}[1]{\left\| #1 \right\|}
\newcommand{\av}[1]{\left| #1 \right|}
\newcommand{\bra}[1]{\left[ #1 \right]}
\newcommand{\set}[1]{\left\{ #1 \right\}}
\newcommand{\dd}{\textnormal{d}}
\newcommand{\sgn}{\textnormal{sgn}}

\newcommand{\cH}{\mathcal{H}}
\newcommand{\cC}{\mathcal{C}}

\newcommand{\bS}{\mathbb{S}}

\newcommand{\cP}{\mathcal{P}}
\newcommand{\bR}{\mathbb{R}}

\def\comm#1#2{{\left\llbracket#1,#2\right\rrbracket}}

\newcommand{\vect}[2]{\pare{\begin{array}{c}#1 \\ #2 \end{array}}}
\newcommand{\ii}{\mathbbm{i}}

\newcommand{\nnorm}[1]{{\left\vert\kern-0.25ex\left\vert\kern-0.25ex\left\vert #1 
    \right\vert\kern-0.25ex\right\vert\kern-0.25ex\right\vert}}

\usepackage[bookmarks,colorlinks]{hyperref}

\newtheorem{theorem}{Theorem}

\theoremstyle{definition}

\title{Global well-posedness and decay for viscous water wave models}

\author[R. Granero-Belinch\'{o}n]{Rafael Granero-Belinch\'{o}n}
\address{Departamento  de  Matem\'aticas,  Estad\'istica  y  Computaci\'on,  Universidad  de Cantabria.  Avda.  Los  Castros  s/n,  Santander,  Spain.}
\email{rafael.granero@unican.es}

\author[S. Scrobogna]{Stefano Scrobogna}
\address{Departamento de An\'alisis Matem\'atico \&  IMUS, Av. de la Reina Mercedes, s/n, 41012 Sevilla, Spain}
\email{scrobogna@us.es}

\begin{document}
\begin{abstract}
The motion of the free surface of an incompressible fluid is a very active research area. Most of these works examine the case of an inviscid fluid. However, in several practical applications, there are instances where the viscous damping needs to be considered. In this paper we derive and study a new asymptotic model for the motion of unidirectional viscous water waves. In particular, we establish the global well-posedness in Sobolev spaces. Furthermore, we also establish the global well-posedness and decay of a fourth order PDE modelling bidirectional water waves with viscosity moving in deep water with or without surface tension effects. \end{abstract}



\maketitle
{\small
\tableofcontents}

\allowdisplaybreaks
\section{Introduction}
\begin{figure}[h]
	\centering
\begin{tikzpicture}[domain=0:2*pi, scale=1]
    \draw (1,3) node { $\Gamma(t)$}; 
    \draw (-0.5,0.5) node { $\Omega(t)$}; 
    \draw[very thick, smooth, variable=\x, blue] plot (\x,{cos(\x r)+2}); 
	\filldraw[blue!10] plot[domain=0:2*pi] (\x,0) -- plot[domain=2*pi:0] (\x,{cos(\x r)+2});
\end{tikzpicture} 
	\caption{Scheme of the problem}\label{Fig1}
\end{figure}
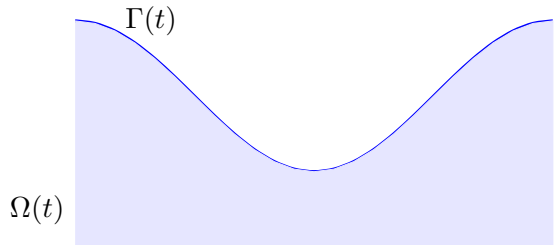

The motion of waves in fluids has been a hot research topic since the XVIIIth century with the works of Laplace and Lagrange. On the one hand there is a large number of papers dealing with the free boundary Euler and Navier-Stokes equations \cite{Beale1981, Lannes2005}. These are free boundary problems and as a consequence the domain of definition $\Omega(t)$ of the functions (the bulk of the fluid) is an unknown of the system that has to be determined from the dynamics (see Figure \ref{Fig1}).

On the other hand the literature on asymptotic models of such free boundary problems is even larger (cf. \cite{Lannes13}). These asymptotic models allow to have very good approximate description of the actual dynamics while simplifying the equations under study. In this direction there are many papers dealing with the case of inviscid fluids and, in particular dealing with asymptotic models for shallow water waves (see for instance \cite{Lannes13} and the references therein) and models of water waves with small steepness (we refer to \cite{matsuno1992nonlinear,matsuno1993two,matsuno1993nonlinear,AkMi2010,aurther2019rigorous} for example). Similar small steepness asymptotics models can be derived for other free boundary problems, such as the Muskat problem \cite{CCGS2013, GGS2020}, see \cite{granero2019asymptotic, GBS2020, Scrobogna2020}. Many of such asymptotic models are used in different applications in Coastal Engineering and Physics. 

Although it is a classical topic, the works studying the case of a viscous fluid are more scarce. The first works studying the case of a viscous water wave date back to Boussinesq \cite{boussinesq1895lois}, Basset \cite{basset1888treatise} and Lamb \cite{lamb1932hydrodynamics}. Since then there are many other papers studying damped water waves. For instance, we refer to the manuscripts of Kakutani \& Matsuuchi \cite{kakutani1975effect}, Ruvinsky \& Freidman \cite{ruvinsky1987fine}, Longuet-Higgins \cite{longuet1992theory}, Jiang, Ting, Perlin \& Schultz \cite{jiang1996moderate}, Joseph \& Wang \cite{joseph2004dissipation}, Wang \& Joseph \cite{wang2006purely} and Wu, Liu \& Yue \cite{wu2006note}. 

According to the work by Dias, Dyachenko \& Zakharov \cite{dias2008theory}, the viscous damping of gravity water waves can be described by the following free boundary problem:
\begin{subequations}\label{eq:DDZ}
\begin{align}
\Delta \phi&=0&&\text{ in }\Omega(t),\\
\rho\left(\phi_t+\frac{1}{2}|\nabla\phi|^2+ Gh\right)&=-2\mu \partial_2^2\phi&&\text{ on }\Gamma(t),\\
h_t&=\nabla\phi\cdot \left(-\partial_{1}h,1\right)+2\frac{\mu}{\rho}\partial_{1}^2h &&\text{ on }\Gamma(t),
\end{align}
\end{subequations}
where $h$ denotes the height of the wave, $\phi$ is the velocity potential and $G,\rho$ and $\mu$ are the gravity acceleration, density and viscosity of the fluid.

Since its appearance, this system was considered by several other authors (see \cite{dutykh2007viscous,dutykh2009visco, dutykh2007dissipative,dutykhA}). The need for simplified asymptotic models for damped water-waves systems was highlighted at first by Longuet-Higgins, which in \cite{longuet1992theory} stated that
\begin{quote}
\emph{For certain applications, however, viscous damping of the waves is important, and it would be highly convenient to have equations and boundary conditions of comparable simplicity as for undamped waves.}
\end{quote}
In this spirit, Kakleas \& Nicholls \cite{kakleas2010numerical} derived a quadratic asymtotic model while Bae, Lin \& Shin \cite{BaeLinShin} derived a cubic asymptotic model  of \eqref{eq:DDZ}. The well-posedness of this quadratic model was studied by Ambrose, Bona \& Nicholls \cite{ambrose2012well} while the well-posedness of the full Dias-Dyachenko-Zakharov was proved by Ngom \& Nicholls in \cite{ngom2018well} in the case of a nonzero surface tension and by Granero-Belinch\'on \& Scrobogna \cite{graneroscrobo3} in the case in which the surface tension can be zero.

In a series of works \cite{graneroscrobo,graneroscrobo2}, the authors, starting with the Dias-Dyachenko-Zakharov \eqref{eq:DDZ} free boundary problem, derived and studied the following bidirectional models of viscous water waves
\begin{equation}\label{NEW2}
\system{\begin{aligned}
& \begin{multlined}
f_{tt}+2\delta\Lambda^2 f_t+ \Lambda f+\beta\Lambda^3 f+\delta^2\Lambda^4 f= \varepsilon\bigg{\lbrace}-\Lambda\left(\left(\cH f_t\right)^2\right)\\
+\partial_x\comm{\cH}{f}\Lambda f 
+\beta\partial_x\comm{\cH}{f}\Lambda^3 f
+\delta\partial_x\comm{\cH}{\cH f _t}\cH \partial_x ^2 f \\
+\delta\Lambda\left(\cH f _t\cH \partial_x ^2 f \right)
-\delta\partial_x\comm{\partial_x^2}{f} \cH f _t +\delta^2\partial_x\comm{\partial_x^2}{f}\Lambda\partial_{x} f\\
-\delta^2\partial_x\comm{\cH}{\partial_{x}^2 f}\partial_{x}^2 f \bigg{\rbrace}, 
\end{multlined}
\\
&f(x,0)=f_0(x),\\
&f_t(x,0)=f_1(x),
\end{aligned}}
\end{equation}
and
\begin{equation}\label{NEW}
\system{\begin{aligned}
& \begin{multlined}
f_{tt}+2\delta\Lambda^2 f_t+ \Lambda f+\beta\Lambda^3 f+\delta^2\Lambda^4 f= \varepsilon\bigg{\lbrace}-\Lambda\left(\left(\cH f_t\right)^2\right)\\
+\partial_x\comm{\cH}{f}\Lambda f 
+\beta\partial_x\comm{\cH}{f}\Lambda^3 f
+\delta\partial_x\comm{\cH}{\cH f _t}\cH \partial_x ^2 f \\
+\delta\Lambda\left(\cH f _t\cH \partial_x ^2 f \right)
-\delta\partial_x\comm{\partial_x^2}{f} \cH f _t\bigg{\rbrace} , 
\end{multlined}
\\
&f(x,0)=f_0(x),\\
&f_t(x,0)=f_1(x),
\end{aligned}}
\end{equation}
where $\varepsilon$ is the steepness parameters that measures the ratio between the anplitude and the wavelength, $\delta>0$ is a dimensionless parameter reflecting the viscous effects, $\beta\geq0$ is the Bond number measuring the ratio between capillary and gravity forces. The operators $\mathcal{H}$ and $\Lambda$ denote the Hilbert transform and the square root of the Laplacian
\begin{align}\label{Hilbert}
\widehat{\mathcal{H}f}(k)=-i\text{sgn}(k) \hat{f}(k) \,, \ \ 
\widehat{\Lambda f}(k)=|k|\hat{f}(k)\,,
\end{align}
and 
$$
\comm{A}{B}f=A(Bf)-B(Af),
$$ is the commutator between two operators acting on the function $f$. In what follows we consider $(x,t)\in\mathbb{S}^1\times[0,T]$ where $\mathbb{S}^1$ denotes the interval $[-\pi,\pi]$ with periodic boundary conditions. Furthermore, we will consider zero-mean initial data $f_0$ and $f_1$.

The purpose of this work is twofold. First we prove that the system \eqref{NEW} is globally well posed for initial data which are sufficiently small. \\

Second we  derive a new asymptotic model of unidirectional viscous water waves. In particular, we  obtain the following nonlocal and nonlinear equation
\begin{multline}\label{NEW4}
2\varepsilon u_{t}=\mathcal{N}u_{x}+2\delta \mathcal{N}u_{xx}+ \mathcal{N}\cH u-\beta \mathcal{N}\cH \partial_x^2 u+\delta^2\mathcal{N}\partial_x^3 u\\
-\varepsilon \mathcal{N}\bigg{\lbrace}2uu_x+\Lambda\comm{\cH}{ \Lambda^{-1} u }u +\beta\Lambda\comm{\cH}{\Lambda^{-1} u }\Lambda^2 u
\\
-\delta\Lambda\comm{\cH}{u}u_x+\delta\partial_x\left(uu_x \right)
+\delta\Lambda\comm{\partial_x^2}{\Lambda^{-1} u} u
\bigg{\rbrace}.
\end{multline}
where the operator 
$$
\mathcal{N}=(1-\delta^2\partial_x^2)^{-1}(1-\delta\partial_x)
$$
is defined in Fourier variables as
$$
\widehat{\mathcal{N}}=\frac{1-\delta ik}{1+\delta^2|k|^2}.
$$

\subsection{Main results}
We start this section introducing some notation that we will use along the paper. We denote with $ C $ any positive constant independent of any physical parameter of the problem. The explicit value of $ C $ may vary from line to line. 

We recall the definition of the homogeneous Sobolev spaces of fractional order 
\begin{align*}
\dot{H}^s=\dot{H}^s\pare{\bS^1} = \set{f\in L^1 \ \left| \ \Lambda^s f \in L^2 \right. }, 
\end{align*}
for any $ s \in \bR $. It is well known that for zero mean function we have that $ H^s =\dot{H}^s $. As both equations preserve the zero mean property from now on we will always use the non-homogeneous notation in order to indicate a Sobolev space of regularity $ s $. Similarly, we define the homogeneous Wiener spaces
\begin{align*}
\dot{A}^s= \dot{A}^s\pare{\bS^1} = \set{f\in L^1 \ \left| \ \widehat{\Lambda^s f} \in \ell^1 \right. }, 
\end{align*}
where $\hat{f}$ denotes the Fourier series of $f$. \\

The first main result of this work is the following theorem:
\begin{theorem}\label{theorem} 
Let $\delta>0$ and $\beta\geq0$. There exists a $ c_0 > 0 $ such that for any $(f_0,f_1)\in H^{6}\times H^{4}$ such that
\begin{equation*}
\norm{f_0}_{H^6} + \norm{f_1}_{H^4} \leq c_0, 
\end{equation*}
 then, there exist a unique global solution $ (f,f_t) $ of \eqref{NEW} stemming from the initial data $ \pare{f_0, f_1} $ which belongs to the energy space
\begin{align*}
f & \in C\left( \bR_+ ;  H^{6}\right), \\
f_t & \in C\left(\bR_+ ;  H^{4}\right)\cap L^2(\bR_+ ; H^{5}).
\end{align*}
Furthermore,
\begin{align*}
\norm{f}_{A^0}+\norm{f_t}_{A^0}&\leq C e^{-t\delta},\\
\norm{f}_{H^r}+\norm{f_t}_{H^s}&\leq C e^{-C(\delta,r,s)t}, &&\forall \  \pare{r, s}\in \left[0, 6\right) \times \left[ 0, 4\right).
\end{align*}
\end{theorem}
Once the local existence and uniqueness was obtained in \cite{graneroscrobo2}, we only need to provide with appropriate energy estimates. In order to do that we will consider a space of low regularity $ X $ and a space of high regularity $ Y $, which will be explicitly defined below. Next we are going to define an energy having the  form
$$
\nnorm{(f,f_t)}_T= \sup _{t\in\bra{0, T}}\set{e^{\alpha t}\norm{\pare{ f\pare{\bullet,  t},f_t \pare{\bullet,  t} } }_{X}}+\|(f,f_t)\|_Y\quad \text{for }\alpha>0.
$$
Equipped with this definition of energy, the rest of the proof is focused on obtaining an inequality of the form
$$
\nnorm{(f,f_t)}_T\leq \cC_0\pare{f_0, f_1}+ P \pare{ \nnorm{(f,f_t)}_T } ,
$$
for certain polynomial $ P $ and constant $\cC_0\pare{f_0, f_1}$ that depends on the initial data. The previous inequality implies that for small enough $\cC_0\pare{f_0, f_1}$, the solution satisfies
$$
\nnorm{\pare{f_0, f_t}}_T\leq 2 \cC_0\pare{f_0, f_1},
$$
for all $ T > 0 $,  then a standard continuation argument allow us to extend the solution to  arbitrary long time intervals. \\

Next, we derive a new asymptotic model of unidirectional viscous water waves. This new model takes the form \eqref{NEW4}. Our second main result is 
\begin{theorem}\label{theorem2} 
Let $\delta>0$ and $\beta\geq0$. Then given and arbitrary zero mean $u_0\in H^{2}$, there exists a unique local strong solution to \eqref{NEW4}
$$
u \in C([0,T^*],H^2)\cap L^2(\bra{0,T^*};H^3),
$$
for a small enough $T^*$ depending only on $\|u_0\|_{H^2}$ and the physical parameters of the problem. Furthermore, there exists a $ c_0 > 0 $ such that for any $u_0\in H^{2}$ satisfying
\begin{equation*}
\norm{u_0}_{H^2} \leq c_0, 
\end{equation*}
then, there exist a unique global solution $u$ of \eqref{NEW4} stemming from the initial data $ u_0 $ which belongs to the energy space
\begin{align*}
u & \in C\left(\bR_+ ;  H^{2}\right)\cap L^2(\bR_+ ; H^{3}).
\end{align*}
Moreover,
\begin{align*}
\norm{u}_{A^0}&\leq C e^{-\frac{\delta}{2} t} ,\\
\norm{u}_{H^r}&\leq C e^{-C(\delta,r)t}, &&\forall \  r\in \left[0, 2\right).
\end{align*}
\end{theorem}
In order to prove the local existence part of this theorem we  use Picard's theorem together with energy estimates in $H^2$ and the commutator structure of part of the nonlinearity. Once the local existence and uniqueness has been obtained, to ensure the global existence and decay we only need to provide with appropriate energy estimates. To do that we are going to define a modified energy $
\nnorm{f}_T$ that has two different contributions. On the one hand we consider the low regularity space $X$ where the solution will decay while on the other hand we will also regard a high regularity space $Y$ where the solution will only remain bounded. The particular choice of $X$ and $Y$ will be clear below. Then the energy will take the form
$$
\nnorm{f}_T= \sup _{t\in\bra{0, T}}\set{e^{\alpha t}\norm{f\pare{\bullet,  t} }_{X}}+\|f\|_Y\quad \text{for }\alpha>0.
$$
Equipped with this definition of energy, the rest of the proof will be devoted to obtain an inequality of the form
$$
\nnorm{u}_T\leq \cC_0\pare{f_0}+ P (\nnorm{u}_T),
$$
for certain polynomial  $ P $ of degree larger than 1, constant $\cC_0\pare{u_0}$ that depends on the initial data. The previous inequality implies that for small enough $\cC_0\pare{u_0}$, the solution satisfies
$$
\nnorm{u}_T\leq 2 \cC_0\pare{u_0},
$$
for all $ T > 0 $,  then a standard continuation argument allow us to extend the solution to  arbitrary long time intervals.

\section{Proof of Theorem \ref{theorem}}
Without loss of generality, we consider $\varepsilon=1$ in \eqref{NEW}. According to the result in \cite{graneroscrobo2}, there is a local in time solution $(f,f_t)$ for the problem \eqref{NEW}.
Let us define the modified energy
\begin{multline*}
\nnorm{(f,f_t)}_{T} =  \ e^{\delta T}\max_{t'\in\bra{0, T}} \set{ \norm{\left(f\pare{t'},f_t\pare{t'}\right)}_{A^0}}  \\
+ \max_{t'\in\bra{0, T}} \set{  \norm{f_t\pare{t'}}_{H^{4}}+ \norm{f\pare{t'}}_{H^{6}} } .
\end{multline*}
The estimates of \cite{graneroscrobo2} assures us moreover that the solution exists at least in a time interval $[0,T_{\text{max}}]$ where $T_{\text{max}} = T_{\max}\pare{f_0, f_1} $ is the maximal lifespan of the solution. 

\subsection{The linear semigroup}
We consider the linear nonhomogeneous problem
\begin{equation}\label{eq:DWW_asint2}
f_{tt}+2\delta\Lambda^2 f_t+ \Lambda f+\beta\Lambda^3 f+\delta^2\Lambda^4 f =F,
\end{equation}
where $F$ is a zero mean forcing. Let us denote with
\begin{align*}
u(x,t) = \vect{f(x,t)}{f_t(x,t)}, && u_0(x) = \vect{f_0(x)}{f_1(x)}, 
\end{align*}
so that \eqref{eq:DWW_asint2} becomes
\begin{align*}
u_t + \mathcal{L} u = \vect{0}{F}, && \mathcal{L} = \pare{ \begin{array}{cc}
0 & -1 \\ \Lambda +\beta\Lambda^3 +\delta^2\Lambda^4 &  2\delta\Lambda^2
\end{array} }. 
\end{align*}
Applying Duhamel principle we write $ u= u_{\textnormal{L}} +  u_{\textnormal{NL}} $ where
\begin{align*}
 u_{\textnormal{L}} \pare{t} = e^{-t\mathcal{L}}u_0, &&  u_{\textnormal{NL}}\pare{t} = \int_0^t e^{-\pare{t-t'} \mathcal{L} }  \vect{0}{F\pare{t'}}\dd t'.
\end{align*}

The eigenvalues of $ L $ are the Fourier multipliers
\begin{equation*}
\lambda_\pm\pare{n} = \delta \av{n}^2 \pm \ii \sqrt{\av{n}\pare{1+\beta\av{n}^2}}, 
\end{equation*}
so we see that the linear operator $ L $ induces both parabolic smoothing effects and oscillating behavior of the solution. Since the solution has zero mean, we have that $\lambda_\pm(n)\neq0$. The two ortonormal eigenvectors associated to $ \lambda_{\pm}\pare{n} $ are
\begin{equation*}
\mathsf{e}_{\pm}\pare{n} = \frac{1}{\sqrt{1+\av{\lambda_{\pm}\pare{n}}^2}} \vect{1}{- \lambda_\pm\pare{n}}, 
\end{equation*}
so that, if we denote
\begin{align*}
D & = \pare{ \begin{array}{cc}
\lambda_- & 0 \\  0  & \lambda_+
\end{array} }, \\
S & =   \pare{ \begin{array}{cc}
1 & 1 \\  - \lambda_-  & - \lambda_+
\end{array} } , \\
S^{-1} & = \frac{1}{\lambda_- - \lambda_+}  \pare{ \begin{array}{cc}
-\lambda_+ & -1 \\   \lambda_- & 1 
\end{array} },
\end{align*}
we have that
$$
e^{-t\mathcal{L}} = S^{-1} e^{-tD} S.
$$
With the above considerations we write $ u_{\textnormal{L}} $ and $ u_{\textnormal{NL}} $ in terms of $ f_0, f_1 $ and $ F $ as
\begin{align*}
\hat{u}_{\textnormal{L}}\pare{t} & = 
\frac{1}{\lambda_- - \lambda_+} 
\vect{\lambda_- e^{-t\lambda_+} - \lambda_+ e^{-t\lambda_-} }{\lambda_- \pare{ e^{-t\lambda_-} -e^{-t\lambda_+} }} \hat{ f }_0\\
&\quad+ \frac{1}{\lambda_- - \lambda_+} \vect{\lambda_+\pare{ e^{-t\lambda_+}-e^{-t\lambda_-} }}{\lambda_- e^{-t\lambda_-} - \lambda_+ e^{-t\lambda_+} } \hat{f}_1
 , \\
\hat{u}_{\textnormal{NL}}\pare{t} & = \int_0^t \frac{1}{\lambda_- - \lambda_+} \vect{\lambda_+\pare{ e^{-\pare{ t-t' }\lambda_+}-e^{-\pare{ t-t' }\lambda_-} }}{\lambda_- e^{-\pare{ t-t' }\lambda_-} - \lambda_+ e^{-\pare{ t-t' }\lambda_+} } \hat{F}\pare{t'}\dd t'. 
\end{align*}

We want to obtain now the decay rates of the linear semigroup. Let us at first check the time-decay of $ u_{\textnormal{L}} $. We can compute that
\begin{equation}\label{eq:estimate_lambda_ratio}
\av{\frac{\lambda_{\pm}}{\lambda_- - \lambda_+}} \leq \frac{1}{2} \pare{  1 +\delta\sqrt{\frac{\av{n}}{1+\beta}} } , 
\end{equation}
since $ \av{n}\geq 1 $ due to conservation of average.
We deduce that, for $ j=0, 1 $
\begin{equation*}
\av{\frac{\lambda_{\pm}\pare{n}}{\lambda_-\pare{n} -\lambda_+\pare{n}} e^{-t\lambda_{\pm}\pare{n}} \hat{f}_j\pare{n}} \leq \frac{ e^{-\delta t} }{2}\pare{1 +\delta\sqrt{\frac{\av{n}}{1+\beta}}} \av{\hat{f}_j\pare{n}}, 
\end{equation*}
which in turn implies that
\begin{equation} \label{eq:control_uL}
\ e^{\delta T}\max_{t'\in\bra{0, T}} \norm{u_{\textnormal{L}}(t')}_{A^0} \leq C\norm{\pare{f_0, f_1}}_{A^{1/2}}. 
\end{equation}
Equivalently, we have that
\begin{equation}
\label{linear_operator_continuity}
\|e^{-t \mathcal{L}}\|_{A^{1/2}\mapsto A^0}\leq C e^{-\delta t}.
\end{equation}

\subsection{Decay in the low regularity space}If we write the equation in its mild formulation using Duhamel's principle, we have that the nonlinear forcing is given by
\begin{equation*}
F= \sum_{j=1}^6 F_j, 
\end{equation*}
where
\begin{align*}
F_1 & = -\Lambda\left(\left(\cH f_t\right)^2\right),\\
F_2 & = \partial_x\comm{\cH}{f}\Lambda f ,\\
F_3 & =  \beta\partial_x\comm{\cH}{f}\Lambda^3 f , \\
F_4 & = \delta\partial_x\comm{\cH}{\cH f _t}\cH \partial_x ^2 f , \\
F_5 & = \delta\Lambda\left(\cH f _t\cH \partial_x ^2 f \right) , \\
F_6 & = \delta\partial_x\comm{\partial_x^2}{f} \cH f _t .
\end{align*}
The goal of the present computations is to provide a control of the form
\begin{align*}
\norm{F\pare{t}}_{A^{1/2}} \leq C e^{-\delta t \pare{1+q}} \nnorm{\pare{f, f_t}}_T^2, && t\in\bra{0, T}, \  q>0. 
\end{align*}

We are going to use the Sobolev embedding 
$$
\|a\|_{A^s}\leq C_\delta\|a\|_{H^{s+1/2+\delta}}\leq C\|a\|_{H^{s+1}},
$$
together with interpolation between Sobolev spaces and the fractional product rule
$$
\|ab\|_{A^s}\leq C_s(\|a\|_{A^0}\|b\|_{A^s}+\|a\|_{A^s}\|b\|_{A^0})\leq C_s\|a\|_{A^s}\|b\|_{A^s},
$$
to estimate $F_j$. We compute
\begin{align*}
\|F_1\|_{A^{1/2}} & \leq C\|f_t\|_{A^0}\|f_t\|_{A^{3/2}}\\
&\leq C\|f_t\|_{A^0}\|f_t\|_{H^{5/2}}\\
&\leq C\|f_t\|_{A^0}^{1+3/8}\|f_t\|_{H^{4}}^{5/8}.
\end{align*}
Using linear interpolation in Wiener spaces
$$
\|a\|_{A^s}\leq C\|a\|_{A^{r}}^{s/r}\|a\|_{A^{0}}^{1-s/r},
$$
we find that
\begin{align*}
\|F_2\|_{A^{1/2}} & \leq \|\comm{\cH}{f}\Lambda f\|_{A^{3/2}}\\
&\leq C(\|f\|_{A^{3/2}}\|f\|_{A^1}+\|f\|_{A^{5/2}}\|f\|_{A^0})\\
&\leq C\|f\|_{A^{5/2}}\|f\|_{A^0}\\
&\leq C\|f\|_{H^{7/2}}\|f\|_{A^0}\\
&\leq C\|f\|_{H^{6}}^{7/12}\|f\|_{A^0}^{1+5/12}.
\end{align*}
Similarly,
\begin{align*}
\|F_3\|_{A^{1/2}} & \leq  C\|\comm{\cH}{f}\Lambda^3 f\|_{A^{3/2}} \\
&\leq C(\|f\|_{A^{3/2}}\|f\|_{A^3}+\|f\|_{A^{3+3/2}}\|f\|_{A^0})\\
&\leq C\|f\|_{A^{9/2}}\|f\|_{A^0}\\
&\leq C\|f\|_{H^{11/2}}\|f\|_{A^0}\\
&\leq C\|f\|_{H^{6}}^{11/12}\|f\|_{A^0}^{1+1/12}.
\end{align*}
Similarly, we have that
\begin{align*}
\|F_4\|_{A^{1/2}} & \leq C\|\comm{\cH}{\cH f _t}\cH \partial_x ^2 f\|_{A^{3/2}}\\
&\leq C\left(\|f_t\|_{A^{0}}\|f\|_{A^{7/2}}+\|f_t\|_{A^{3/2}}\|f\|_{A^{2}}\right)\\
&\leq C\left(\|f_t\|_{A^{0}}\|f\|_{H^{6}}+\|f_t\|_{A^{0}}^{1/2}\|f_t\|_{A^{3}}^{1/2}\|f\|_{A^{0}}^{3/5}\|f\|_{A^{5}}^{2/5}\right)\\
&\leq  C\left(\|f_t\|_{A^{0}}\|f\|_{H^{6}}+\|f_t\|_{A^{0}}^{1/2}\|f_t\|_{H^{4}}^{1/2}\|f\|_{A^{0}}^{3/5}\|f\|_{H^{6}}^{2/5}\right) ,
\end{align*}
\begin{align*}
\|F_5\|_{A^{1/2}} & \leq C\left(\|f_t\|_{A^{0}}\|f\|_{H^{6}}+\|f_t\|_{A^{0}}^{1/2}\|f_t\|_{H^{4}}^{1/2}\|f\|_{A^{0}}^{3/5}\|f\|_{H^{6}}^{2/5}\right) , \\
\end{align*}
Using the commutator structure together with the product rule in Wiener spaces, we estimate
\begin{align*}
\|F_6\|_{A^{1/2}} & \leq C\|\partial_x^3f \cH f _t+3\partial_x^2f \partial_x\cH f _t+\partial_xf \partial_x^2\cH f _t\|_{A^{1/2}}\\
&\leq C\bigg{(}\|f\|_{A^{3+1/2}}\|f _t\|_{A^0}+\|f\|_{A^{3}}\|f _t\|_{A^{1/2}}\\
&\quad+\|f\|_{A^{2+1/2}}\|f _t\|_{A^1}+\|f\|_{A^{2}}\|f _t\|_{A^{1+1/2}}\\
&\quad+\|f\|_{A^{1+1/2}}\|f _t\|_{A^2}+\|f\|_{A^{1}}\|f _t\|_{A^{2+1/2}}\bigg{)}.
\end{align*}
Using interpolation in Wiener spaces and then the Sobolev embedding
\begin{align*}
\|f\|_{A^{5+2/5}}\leq C\|f\|_{H^{6}} && \text{ and } && \|f_t\|_{A^{3+2/5}}\leq C\|f_t\|_{H^{4}} , 
\end{align*}
we compute that
\begin{align*}
\|F_6\|_{A^{1/2}} 
\leq & \  C\bigg{(}\|f\|_{A^0}^{19/54}\|f\|_{H^{6}}^{35/54}\|f _t\|_{A^0}+\|f\|_{A^0}^{4/9}\|f\|_{H^{6}}^{5/9}\|f _t\|_{A^0}^{29/34}\|f _t\|_{H^{4}}^{5/34}\\
&\ +\|f\|_{A^0}^{29/54}\|f\|_{H^{6}}^{25/54}\|f _t\|_{A^0}^{12/17}\|f _t\|_{H^{4}}^{5/17}\\
&\ +\|f\|_{A^0}^{17/27}\|f\|_{H^{6}}^{10/27}\|f _t\|_{A^0}^{19/34}\|f _t\|_{H^{4}}^{15/34}\\
&\ +\|f\|_{A^0}^{13/18}\|f\|_{H^{6}}^{5/18}\|f _t\|_{A^0}^{7/10}\|f _t\|_{H^{4}}^{10/17}\\
&\ +\|f\|_{A^0}^{22/27}\|f\|_{H^{6}}^{5/27}\|f _t\|_{A^0}^{9/34}\|f _t\|_{H^{4}}^{25/34}\bigg{)}.
\end{align*}

Let us recall that using Duhamel formulation the solution then can be written as
$$
u(x,t)=e^{-tL}u_0+\int_0^t e^{-\pare{t-t'} L }  \vect{0}{F\pare{t'}}\dd t'
$$
and satifies, 
$$
\av{ \hat{u}(n,t)}\leq Ce^{-\delta t}(1+\sqrt{|n|}) \av{\hat{f}_j(n)}+C \int_0^t e^{-\pare{t-t'}\delta n^2} \sqrt{\av{n}} \av{\hat{F}\pare{n, t'}} \dd t'. 
$$
Using
$$
-\pare{t-t'}\delta n^2\leq -\pare{t-t'}\delta \leq 0,
$$
we can estimate
$$
\|u(t)\|_{A^{0}}\leq Ce^{-\delta t}\|(f_0,f_1)\|_{A^{1/2}}+C e^{-\delta t}\int_0^t e^{\delta t'}\sum_{j=1}^6\|F_j(t')\|_{A^{1/2}}\dd t'.
$$
Recalling the previous estimates for $\|F_j(t')\|_{A^{1/2}}$ and the definition of the norm $\nnorm{(f,f_t)}_T$, we have that
\begin{align*}
\norm{F \pare{t'}}_{A^{1/2}}\leq C e^{-\delta t'(1+73/918)}\nnorm{(f,f_t)}_T^2, && 0\leq t' \leq t \leq T.
\end{align*}
We conclude that
\begin{multline}
\label{eq:loworder}
 e^{\delta t}\max_{t'\in\bra{0, t}} \set{ \norm{\left(f\pare{t'},f_t\pare{t'}\right)}_{A^0}} \\
 \begin{aligned}
\leq & \  C\|(f_0,f_1)\|_{A^{1/2}}+C\nnorm{(f,f_t)}_T^2\int_0^t e^{-(73/918)\delta t'}\dd t'\\
 \leq & \  C\|(f_0,f_1)\|_{A^{1/2}}+C\nnorm{(f,f_t)}_T^2
\end{aligned}
\end{multline}

\subsection{Boundedness in the high regularity space}
Similarly as in \cite{graneroscrobo2}, we test the equation against $ \Lambda^8 f_t $, integrate in $ \bS^1 $ and integrate by parts obtaining the energy balance 
 \begin{equation}\label{eq:energy_balance1}
 \frac{1}{2}\frac{\dd}{\dd t}\mathfrak{E}(t)+\mathfrak{D}\pare{t'}=\sum_{i=1}^6 I_i\pare{t},
 \end{equation}
with 
\begin{equation*}
\begin{aligned}
\mathfrak{E} (t) & =\norm{f_t\pare{t'}}_{H^{4}}^2+\beta \norm{f\pare{t'}}_{H^{4+3/2}}^2+\delta^2 \norm{f\pare{t'}}_{H^{6}}^2+ \norm{f\pare{t'}}_{H^{4+1/2}}^2, \\
\mathfrak{D} (t) & =2\delta\norm{f_t(t)}_{H^{5}}^2,
\end{aligned}
\end{equation*}
and
\begin{align*}
I_1 \pare{t} &=-\int_{\mathbb{S}^1}\Lambda\left(\left(\cH f_t\right)^2\right)\Lambda^7 f_t \ \text{d} x \\
I_2 \pare{t} &=\int_{\mathbb{S}^1}\partial_x\comm{\cH}{f}\Lambda f\Lambda^8 f_t \ \text{d} x  \\
I_3 \pare{t} &=\beta\int_{\mathbb{S}^1}\partial_x\comm{\cH}{f}\Lambda^3 f\Lambda^8 f_t \ \text{d} x \\
I_4 \pare{t} &=\delta\int_{\mathbb{S}^1}\partial_x\comm{\cH}{\cH f _t}\cH \partial_x ^2 f\Lambda^8 f_t \ \text{d} x  \\
I_5 \pare{t} &=\delta\int_{\mathbb{S}^1}\Lambda\left(\cH f _t\cH \partial_x ^2 f \right)\Lambda^8 f_t \ \text{d} x \\
I_6 \pare{t} &=-\delta\int_{\mathbb{S}^1}\partial_x\comm{\partial_x^2}{f} \cH f _t\Lambda^8 f_t \ \text{d} x .
\end{align*}

Using the self-adjointness of the operator $\Lambda$ together with H\"older's inequality and the Sobolev embedding
$$
\|g\|_{L^4}\leq C\|g\|_{H^{0.25}},
$$
we find that
\begin{align}
I_1(t)&=-\int_{\mathbb{S}^1}\left(\left(\cH f_t\right)^2\right)\Lambda^9 f_t \ \text{d} x \nonumber\\
&=-\int_{\mathbb{S}^1}\left(\left(\cH f_t\right)^2\right)\partial_x^4 \Lambda^5 f_t \ \text{d} x \nonumber\\
&=-\int_{\mathbb{S}^1}\partial_x^4\left(\left(\cH f_t\right)^2\right)\Lambda^5 f_t \ \text{d} x\nonumber\\
&=-\int_{\mathbb{S}^1}\left(2\cH f_t\Lambda\partial_x^3 f_t+6(\Lambda\partial_x f_t)^2+8\Lambda f_t\partial_x^2\Lambda f_t\right)\Lambda^5 f_t \ \text{d} x\nonumber\\
&\leq C\|f_t\|_{H^5}\left(\|f_t\|_{H^4}\|\cH f_t\|_{L^\infty}+\|f_t\|_{H^{2.25}}^2+\|f_t\|_{H^3}\|\Lambda f_t\|_{L^\infty}\right)\nonumber\\
&\leq C\|f_t\|_{H^{5}}\|f_t\|_{H^{4}}\|f_t\|_{H^{2.25}}\nonumber\\
&\leq \sigma\|f_t\|_{H^{5}}^2+C\|f_t\|_{H^{4}}^2\|f_t\|_{H^{2.25}}^2\ ,\label{ineq:I1}
\end{align}
for $\sigma>0$ to be fixed below.

Furthermore, using interpolation between Sobolev spaces, the embedding
$$
H^s\subset H^{r}\,,r\leq s, 
$$
we obtain the estimate
\begin{equation}\label{eq:I1}
I_1\pare{t}\leq \sigma\|f_t\|_{H^{5}}^2+C\nnorm{(f,f_t)}_T^4 e^{-(\delta/2) t}. 
\end{equation}

Equipped with \eqref{commutatorH}, we can estimate $I_2$ as follows
\begin{align}
I_2(t)= & \ \int_{\mathbb{S}^1}\Lambda^4\partial_x\comm{\cH}{f}\Lambda f\Lambda^4 f_t \ \text{d} x \  \nonumber\\
\leq &\ \norm{ \partial_x^5\comm{\cH}{f}\Lambda f}_{L^2}\norm{ \Lambda^4 f_t}_{L^2} \nonumber\\
\leq & \ \norm{\partial_x^5f}_{L^\infty}\norm{ \Lambda f}_{L^2}\norm{ \Lambda^4 f_t}_{L^2} \nonumber. 
\end{align}
As a consequence, by interpolation in Wiener and Sobolev spaces, we have that
\begin{align}\label{eq:I2}
I_2\pare{t}\leq & \ C\nnorm{(f,f_t)}_T^3 e^{-(\delta/2) t} .
\end{align}
Analogously, we find that
\begin{align} \label{eq:I3}
I_3\pare{t}&\leq \norm{\partial_x^4f}_{L^\infty}\norm{ \Lambda^3 f}_{L^2}\norm{ \Lambda^5 f_t}_{L^2}\nonumber\\
&\leq \norm{\partial_x^4f}_{L^\infty}^2\norm{ \Lambda^3 f}_{L^2}^2+\sigma\norm{ \Lambda^5 f_t}_{L^2}^2\nonumber\\
&\leq  C\nnorm{(f,f_t)}_T^4 e^{-(\delta/2) t}+\sigma\norm{ \Lambda^5 f_t}_{L^2}^2.
\end{align}
 
We can decompose $I_4(t)$ as follows
\begin{align*}
I_4(t)&=\delta\int_{\mathbb{S}^1}\left[\Lambda(\cH f_t\Lambda\partial_x f)+\partial_x(\cH f_t\partial_x^2 f)\right]\Lambda^{8} f_t \ \text{d} x \\
&=\delta\int_{\mathbb{S}^1}\left[\Lambda(\cH f_t\Lambda\partial_x f)-\partial_x(\cH f_t\Lambda^2 f)\right]\Lambda^{8} f_t \ \text{d} x \  \\
&=J_1^4+J_2^4,
\end{align*}
with
\begin{align*}
J_1^4&=\delta\int_0^t\int_{\mathbb{S}^1}\Lambda(\cH f_t\Lambda\partial_x f)\Lambda^{8} f_t \ \text{d} x \  \dd t'\\
J_2^4&=-\delta\int_0^t\int_{\mathbb{S}^1}\partial_x(\cH f_t\Lambda^2 f)\Lambda^{8} f_t \ \text{d} x \  \dd t'.
\end{align*}
We will use the fractional Leibniz rule (see \cite{grafakos2014kato,kato1988commutator,kenig1993well}):
$$
\|\Lambda^s(uv)\|_{L^p}\leq C\left(\|\Lambda^s u\|_{L^{p_1}}\|v\|_{L^{p_2}}+\|\Lambda^s v\|_{L^{p_3}}\|u\|_{L^{p_4}}\right),
$$
which holds whenever
$$
\frac{1}{p}=\frac{1}{p_1}+\frac{1}{p_2}=\frac{1}{p_3}+\frac{1}{p_4}\qquad \mbox{where $1/2<p<\infty,1<p_i\leq\infty$},
$$
and $s>\max\{0,1/p-1\}$. Using the fractional Leibniz rule and the self-adjointness of the operator $\Lambda$, we compute
\begin{align*}
J_1^4\pare{t}&=\delta\int_{\mathbb{S}^1}\Lambda^{4}(\cH f_t\Lambda\partial_x f)\Lambda^{5} f_t \ \text{d} x \nonumber\\
&\leq\delta \|\Lambda^{4}(\cH f_t\Lambda\partial_x f)\|_{L^2}\|\Lambda^{5} f_t\|_{L^{2}}\nonumber\\
&\leq\delta C(\|f_t\|_{H^{1}}\|f\|_{H^{6}}+\|f_t\|_{H^{4}}\|f\|_{H^{3}})\|f_t\|_{H^{5}}\nonumber\\
&\leq\delta C(\|f_t\|_{H^{1}}^2\|f\|_{H^{6}}^2+\|f_t\|_{H^{4}}^2\|f\|_{H^{3}}^2)+\sigma\|f_t\|_{H^{5}}^2\nonumber\\
&\leq C\nnorm{f}_T^4 e^{-(\delta/2)t} +\sigma\|f_t\|_{H^{5}}^2.
\end{align*}
The terms $J^4_2$ and $I_5=J^4_1$ can be estimated in a similar way and we find that
\begin{equation}\label{eq:I4I5}
I_4 \pare{t} + I_5\pare{t} \leq  C\nnorm{f}_T^4 e^{-(\delta/2)t} +\sigma\|f_t\|_{H^{5}}^2.
\end{equation}
Now we are left with $I_6$. We remark that
\begin{equation*}
I_6(t) =-\delta\int_{\mathbb{S}^1}\partial_x\left[\partial_x^2f\cH f _t+2\partial_xf\Lambda f _t\right]\partial_x^4\Lambda^4 f_t \ \text{d} x .
\end{equation*}
Integrating by parts, we find that
\begin{equation*}
I_6(t) =\delta\int_{\mathbb{S}^1}\partial_x^4\left[\partial_x^2f\cH f _t+2\partial_xf\Lambda f _t\right]\partial_x\Lambda^4 f_t \ \text{d} x .
\end{equation*}
Hence, using the same ideas as before, we have that 
\begin{equation*}
I_6(t) \leq \nnorm{(f,f_t)}_T^4e^{-(\delta/2)t}+\sigma\|f_t\|_{H^5}^2+ 2\int_{\mathbb{S}^1}\partial_xf\Lambda^5 f _t\partial_x\Lambda^4 f_t \ \text{d} x .
\end{equation*}
The term
$$
J^6_1=2\int_{\mathbb{S}^1}\partial_xf\Lambda^5 f _t\partial_x\Lambda^4 f_t \ \text{d} x
$$
is the highest order term. However, it has an inner commutator structure that we can exploit as follows:
\begin{align*}
J^6_1(t) &=\int_{\mathbb{S}^1}\cH(\partial_xf\Lambda^5 f _t)\Lambda^5 f_t \ \text{d} x -\int_{\mathbb{S}^1}\partial_xf\Lambda^5 f _t\cH \Lambda^5 f_t \ \text{d} x\\
&=\int_{\mathbb{S}^1}\comm{\cH}{\partial_xf}\Lambda^5 f _t\Lambda^5 f_t \ \text{d} x.
\end{align*}
Then, recalling \eqref{commutatorH}, we conclude that
\begin{equation}\label{eq:I6}
I_6(t) \leq \nnorm{(f,f_t)}_T^4e^{-(\delta/2)t}+\sigma\|f_t\|_{H^5}^2.
\end{equation}

\subsection{Finishing the proof of Theorem \ref{theorem}} \label{sec:absurdum}
Collecting \eqref{eq:I1}, \eqref{eq:I2}, \eqref{eq:I3}, \eqref{eq:I4I5} and \eqref{eq:I6} and taking $\sigma$ small enough, we conclude 
 \begin{equation}\label{eq:energy_balance1v2}
\frac{\dd}{\dd t}\mathfrak{E}(t)+\mathfrak{D}\pare{t'}\leq C\nnorm{(f,f_t)}_T^2e^{-(\delta/2)t}+C\nnorm{(f,f_t)}_T^4e^{-(\delta/2)t}.
 \end{equation}
Integrating in time and using \eqref{eq:loworder}, we conclude the polynomial bound
\begin{multline*}
\nnorm{(f,f_t)}_T+\int_0^T\mathfrak{D}\pare{t'}\ \dd t'\\
\leq C\|(f_0,f_1)\|_{A^{1/2}}+\mathfrak{E}(0)+C\pare{ \nnorm{(f,f_t)}_T^2+\nnorm{(f,f_t)}_T^4 },
\end{multline*}
thus, there exists a (fixed, positive) constant $1<C^+$ such that
\begin{equation}
\label{eq:enineqCast}
\nnorm{(f,f_t)}_T\leq C^+ \bra{ \pare{ \norm{f_0}_{H^6}+ \norm{f_1}_{H^{4}} }+\pare{ \nnorm{(f,f_t)}_T^2+\nnorm{(f,f_t)}_T^4 }}.
\end{equation}
We observe that, in the previous estimates, we have not used any hypothesis on the size of the initial data and the previous bound is valid for every solution and $ T\in\pare{0, T_{\max}} $. 

We want to prove that, there exists a $ c_0 > 0 $ such that for any solution of \eqref{NEW} stemming from an initial data
\begin{equation}\label{eq:smallness_c0_hyp}
 \norm{f_0}_{H^6}+ \norm{f_1}_{H^{4}} \leq c_0, 
\end{equation}
the inequality
\begin{equation*}
\nnorm{\pare{f, f_t}}_T \leq C\pare{ \norm{f_0}_{H^6}+ \norm{f_1}_{H^{4}} }, 
\end{equation*}
holds true for any $ T > 0 $ and thus the solution is global by a standard continuation argument. 

Let us assume that the solution does not stay bounded for all times, the contrary being true would imply that the solution is global by a continuation argument.  If the initial data is small enough, we  can find $T$ such that
$$
\nnorm{\pare{f, f_t}}_T = \frac{3}{4} < 1.
$$
The above inequality allow us to deduce  the polynomial bound
\begin{equation}
\label{eq:enineqCast2}
\nnorm{(f,f_t)}_T\leq 2C^+ \bra{ \pare{ \norm{f_0}_{H^6}+ \norm{f_1}_{H^{4}} }+ \nnorm{(f,f_t)}_T^2},
\end{equation}
or equivalently,
\begin{equation}
\label{eq:enineqCast2v2}
2C^+\nnorm{(f,f_t)}_T\leq (2C^+)^2  \pare{ \norm{f_0}_{H^6}+ \norm{f_1}_{H^{4}} }+ \bra{2C^+\nnorm{(f,f_t)}_T}^2.
\end{equation}
Then, without loss of generality we can restrict our analysis to a polynomial of the form
\begin{equation}
\label{eq:enineqCast3}
\nnorm{(f,f_t)}_T\leq \cC_0\pare{f_0, f_1}+ \nnorm{(f,f_t)}_T^2.
\end{equation}
Now we observe that if 
$$
\mathcal{C}_0\ll 1
$$
is small enough, the polynomial
$$
\mathcal{Q}(y)=\mathcal{C}_0-y+y^2
$$
has two positive real roots 
\begin{equation*}
y_{\pm} = \frac{1\pm \sqrt{1-4\cC_0}}{2}, 
\end{equation*}
moreover if $ 0<\cC_0 \ll 1 $
$$
 y_-=\min\{y_+,y_-\}  =\frac{1 - \sqrt{1-4\cC_0}}{2}\leq 2\mathcal{C}_0.
$$
Furthermore, analogously as in in \cite{graneroscrobo2}, we know that the application $ t\mapsto  \nnorm{\pare{f, f_t}}_t $ is continuous for $ t\in \left[0, T_{\max}\right) $. This, together with the smallness in the initial data, implies that
\begin{equation*}
\nnorm{\pare{f, f_t}}_{T} \in \bra{0, y_-}. 
\end{equation*}
We combine the above deduction with the estimate $ y_- \leq 2\cC_0 $ and we deduce that
$$
\nnorm{(f,f_t)}_{T}\leq 2\cC_0<\frac{3}{4},
$$
if we take $\cC_0$ small enough.  This is a contradiction with the definition of $ T $ and implies that the solution is global.

\section{Derivation of (\ref{NEW4})}
Our starting point in this section is \eqref{NEW2}:
\begin{multline}\label{NEW3}
f_{tt}=-2\delta\Lambda^2 f_t- \Lambda f-\beta\Lambda^3 f-\delta^2\Lambda^4 f\\
+\varepsilon \bigg{\lbrace} -\Lambda\left(\left(\cH f_t\right)^2\right)+\partial_x\comm{\cH}{f}\Lambda f +\beta\partial_x\comm{\cH}{f}\Lambda^3 f
\\
+\delta\partial_x\comm{\cH}{\cH f _t}\cH \partial_x ^2 f+\delta\Lambda\left(\cH f _t\cH \partial_x ^2 f \right)
+\delta^2\partial_x\comm{\partial_x^2}{f}\Lambda\partial_{x} f\\
-\delta\partial_x\comm{\partial_x^2}{f} \cH f _t
-\delta^2\partial_x\comm{\cH}{\partial_{x}^2 f}\partial_{x}^2 f \bigg{\rbrace}.
\end{multline}
Let us introduce the 'far-field' variables,
$$
\chi=x-t,\quad \tau=\varepsilon t.
$$
Then, we have that
$$
\frac{\partial}{\partial t}f(\chi(x,t),\tau(t))=-f_{\chi}+\varepsilon f_{\tau},
$$
and
$$
\frac{\partial^2}{\partial t^2}f(\chi(x,t),\tau(t))=f_{\chi\chi}-2\varepsilon f_{\tau\chi}+\varepsilon^2 f_{\tau\tau}.
$$
After neglecting terms of $O(\varepsilon^2)$, \eqref{NEW3} reads
\begin{multline*}
\left(f_{\chi}-2\varepsilon f_{\tau}\right)_\chi=-2\delta\Lambda^2 (-f_{\chi}+\varepsilon f_{\tau})- \Lambda f-\beta\Lambda^3 f-\delta^2\Lambda^4 f\\
+\varepsilon \bigg{\lbrace}-\Lambda\left(\left(\cH f_\chi\right)^2\right)+\partial_\chi\comm{\cH}{f}\Lambda f +\beta\partial_\chi\comm{\cH}{f}\Lambda^3 f
\\
-\delta\partial_\chi\comm{\cH}{\cH f _\chi}\cH \partial_\chi ^2 f-\delta\Lambda\left(\cH f _\chi\cH \partial_\chi ^2 f \right)
+\delta^2\partial_\chi\comm{\partial_x^2}{f}\Lambda\partial_{\chi} f\\
+\delta\partial_\chi\comm{\partial_x^2}{f} \cH f _\chi
-\delta^2\partial_\chi\comm{\cH}{\partial_{x}^2 f}\partial_{\chi}^2 f\bigg{\rbrace}.
\end{multline*}

Integrating in $\chi$ and using our previous notation for the space and time variables we find the equation
\begin{multline*}
f_{x}-2\varepsilon f_{t}=2\delta\partial_x (-f_{x}+\varepsilon f_{\tau})- \cH f+\beta\cH \partial_x^2 f-\delta^2\partial_x^3 f\\
+\varepsilon \bigg{\lbrace}-\cH\left(\left(\cH f_x\right)^2\right)+\comm{\cH}{f}\Lambda f +\beta\comm{\cH}{f}\Lambda^3 f
\\
-\delta\comm{\cH}{\Lambda f}\cH \partial_x ^2 f-\delta\cH\left(\Lambda f\cH \partial_x ^2 f \right)
+\delta^2\comm{\partial_x^2}{f}\Lambda\partial_{x} f\\
+\delta\comm{\partial_x^2}{f} \cH f _x
-\delta^2\comm{\cH}{\partial_{x}^2 f}\partial_{x}^2 f\bigg{\rbrace}.
\end{multline*}

Regrouping terms we can equivalently write
\begin{multline}\label{eq:aux}
(1+\delta\partial_x)2\varepsilon f_{t}=f_{x}+2\delta f_{xx}+ \cH f-\beta\cH \partial_x^2 f+\delta^2\partial_x^3 f\\
-\varepsilon \bigg{\lbrace}-\cH\left(\left(\cH f_x\right)^2\right)+\comm{\cH}{f}\Lambda f +\beta\comm{\cH}{f}\Lambda^3 f
\\
-\delta\comm{\cH}{\Lambda f}\cH \partial_x ^2 f-\delta\cH\left(\Lambda f\cH \partial_x ^2 f \right)
+\delta^2\comm{\partial_x^2}{f}\Lambda\partial_{x} f\\
+\delta\comm{\partial_x^2}{f} \cH f _x
-\delta^2\comm{\cH}{\partial_{x}^2 f}\partial_{x}^2 f\bigg{\rbrace}.
\end{multline}

We observe that, taking the operator
$$
\mathcal{N}=(1-\delta^2\partial_x^2)^{-1}(1-\delta\partial_x),
$$
we find that
\begin{multline}\label{NEW5}
2\varepsilon f_{t}=\mathcal{N}f_{x}+2\delta \mathcal{N}f_{xx}+ \mathcal{N}\cH f-\beta \mathcal{N}\cH \partial_x^2 f+\delta^2\mathcal{N}\partial_x^3 f\\
-\varepsilon \mathcal{N}\bigg{\lbrace}-\cH\left(\left(\cH f_x\right)^2\right)+\comm{\cH}{f}\Lambda f +\beta\comm{\cH}{f}\Lambda^3 f
\\
-\delta\comm{\cH}{\Lambda f}\cH \partial_x ^2 f-\delta\cH\left(\Lambda f\cH \partial_x ^2 f \right)
+\delta^2\comm{\partial_x^2}{f}\Lambda\partial_{x} f\\
+\delta\comm{\partial_x^2}{f} \cH f _x
-\delta^2\comm{\cH}{\partial_{x}^2 f}\partial_{x}^2 f\bigg{\rbrace}.
\end{multline}
As in \cite{granero2021motion}, we now define 
$$
u=\Lambda f.
$$
This new unknown solves the following equation
\begin{multline*}
2\varepsilon u_{t}=\mathcal{N}u_{x}+2\delta \mathcal{N}u_{xx}+ \mathcal{N}\cH u-\beta \mathcal{N}\cH \partial_x^2 u+\delta^2\mathcal{N}\partial_x^3 u\\
-\varepsilon \Lambda\mathcal{N}\bigg{\lbrace}-\cH\left(u^2\right)+\comm{\cH}{ \Lambda^{-1}u }u +\beta\comm{\cH}{\Lambda^{-1}u}\Lambda^2 u
\\
+\delta\comm{\cH}{u}\cH \Lambda u+\delta\cH\left(u\cH \Lambda u \right)
+\delta^2\comm{\partial_x^2}{\Lambda^{-1}u}u_x\\
+\delta\comm{\partial_x^2}{\Lambda^{-1}u} u
-\delta^2\comm{\cH}{\Lambda u}\Lambda u\bigg{\rbrace}.
\end{multline*}
The previous equation can be written equivalently as 
\begin{multline}\label{NEW6}
2\varepsilon u_{t}=\mathcal{N}u_{x}+2\delta \mathcal{N}u_{xx}+ \mathcal{N}\cH u-\beta \mathcal{N}\cH \partial_x^2 u+\delta^2\mathcal{N}\partial_x^3 u\\
-\varepsilon \mathcal{N}\bigg{\lbrace}2uu_x+\Lambda\comm{\cH}{\Lambda^{-1}u}u +\beta\Lambda\comm{\cH}{\Lambda^{-1}u}\Lambda^2 u
\\
-\delta\Lambda\comm{\cH}{u}u_x+\delta\partial_x\left(uu_x \right)
+\delta^2\Lambda\comm{\partial_x^2}{\Lambda^{-1}u}u_x\\
+\delta\Lambda\comm{\partial_x^2}{\Lambda^{-1}u} u
-\delta^2\Lambda\comm{\cH}{\Lambda u}\Lambda u\bigg{\rbrace}.
\end{multline}

Neglecting now the nonlinear terms that are $O(\varepsilon \delta^2)$ we conclude \eqref{NEW4}.

\section{Proof of Theorem \ref{theorem2}}
\subsection{Local well-posedness}
First, we observe that \eqref{NEW4} can be written as
\begin{multline}\label{NEW4v2}
2\varepsilon u_{t}=\mathcal{N}u_{x}+2\delta \mathcal{N}u_{xx}+ \mathcal{N}\cH u-\beta \mathcal{P}\cH \partial_x^2 u+\beta\delta \mathcal{P}\Lambda \partial_x^2 u+\delta^2\mathcal{P}\partial_x^3 u-\delta^3\mathcal{P}\partial_x^4 u\\
-\varepsilon \mathcal{N}\bigg{\lbrace}2uu_x+\Lambda\comm{\cH}{\Lambda^{-1}u}u +\beta\Lambda\comm{\cH}{\Lambda^{-1}u}\Lambda^2 u
\\
-\delta\Lambda\comm{\cH}{u}u_x+\delta\partial_x\left(uu_x \right)
+\delta\Lambda\comm{\partial_x^2}{\Lambda^{-1}u} u
\bigg{\rbrace},
\end{multline}
where the operator 
$$
\mathcal{P}=(1-\delta^2\partial_x^2)^{-1}
$$
is defined in Fourier variables as
$$
\widehat{\mathcal{P}}=\frac{1}{1+\delta^2|k|^2}.
$$
Then, using
$$
\mathcal{P}=Id+\mathcal{P}\delta^2\partial_x^2
$$
 we can observe that the terms
$$
\beta\delta \mathcal{P}\Lambda \partial_x^2 u=-\frac{\beta}{\delta}\Lambda u+ \frac{\beta}{\delta}\mathcal{P}\Lambda u
$$
and
$$
-\delta^3\mathcal{P}\partial_x^4 u=\delta\partial_x^2 u-\delta\mathcal{P}\partial_x^2 u
$$
are of parabolic type.

To simplify the notation, in the course of this proof we take $\varepsilon=1$. Now we obtain the \emph{a priori} estimates in the $H^2$ Sobolev space. These estimates implies the local existence of solution after a standard regularization approach using the periodic heat kernel as mollifier. 

We start noticing that
$$
\int_{\mathbb{S}^1} u(x,t)\dd x=\int_{\mathbb{S}^1} u(x,0)\dd x=0.
$$

Now we test \eqref{NEW4v2} against $\Lambda^{4}u$. Then we obtain that
\begin{align*}
\frac{\dd}{\dd t}\|u\|_{H^2}^2&=L+NL_1+NL_2+NL_3+NL_4+NL_5+NL_6,
\end{align*}
where
\begin{align*}
L&=\int_{\mathbb{S}^1}\lbrace\mathcal{N}u_{x}+2\delta \mathcal{N}u_{xx}+ \mathcal{N}\cH u-\beta \mathcal{P}\cH \partial_x^2 u\\
&\qquad+\beta\delta \mathcal{P}\Lambda \partial_x^2 u+\delta^2\mathcal{P}\partial_x^3 u-\delta^3\mathcal{P}\partial_x^4 u\rbrace\Lambda^4 u \dd x,\\
NL_1&=-2\int_{\mathbb{S}^1}\mathcal{N}(uu_x)\Lambda^4 u \dd x,\\
NL_2&=-\int_{\mathbb{S}^1}\mathcal{N}\Lambda\comm{\cH}{\Lambda^{-1}u }u \Lambda^4 u \dd x,\\
NL_3&=-\beta\int_{\mathbb{S}^1}\mathcal{N}\Lambda\comm{\cH}{\Lambda^{-1}u}\Lambda^2 u\Lambda^4 u \dd x,\\
NL_4&=\delta\int_{\mathbb{S}^1}\mathcal{N}\Lambda\comm{\cH}{u}\partial_x u \Lambda^4 u \dd x,\\
NL_5&=-\delta\int_{\mathbb{S}^1}\mathcal{N}\partial_x\left(uu_x \right)\Lambda^4 u \dd x,\\
NL_6&=-\delta\int_{\mathbb{S}^1}\mathcal{N}\Lambda\comm{\partial_x^2}{\Lambda^{-1}u} u \Lambda^4 u \dd x.
\end{align*}
After a number of integrations by parts, we find that
\begin{align*}
L&=\int_{\mathbb{S}^1}\left(\delta\mathcal{P}u_{xx}-\delta \mathcal{P}\Lambda u+\beta\delta \mathcal{P}\Lambda u_{xx}-\delta^3\mathcal{P}\partial_x^4 u\right)\Lambda^4 u \dd x\\
&=-\delta \|\mathcal{P}^{1/2}u_{xxx}\|_{L^2}^2-\delta \|\mathcal{P}^{1/2} \Lambda^{5/2}u\|_{L^2}^2-\delta\beta \|\mathcal{P}^{1/2} \Lambda^{1/2}u_{xxx}\|_{L^2}^2-\delta^3 \| \mathcal{P}^{1/2} \Lambda^{4}u\|_{L^2}^2.
\end{align*}
Furthermore, using the parabolic character of some of the terms in $L$, we find that
\begin{align*}
L&\leq -\frac{1}{\delta} \|u_{xx}\|_{L^2}^2-\delta \|  \Lambda^{3}u\|_{L^2}^2+C\|u\|_{H^2}^2.
\end{align*}
For the first nonlinear term $NL_1$, we integrate by parts and use that $\mathcal{N}$ can absorb one derivative to find the estimate
$$
NL_1\leq \|u\|_{H^3}\|u^2\|_{H^1}\leq C \|u\|_{H^3}\|u\|_{H^1}^2\leq \sigma\|u\|_{H^3}^2+C\|u\|_{H^1}^4,
$$
for $\sigma>0$ that will be fixed later.

We recall the following commutator estimate (see equation (1.13) in \cite{dawson2008decay})
\begin{align}\label{commutatorH}
\norm{ \partial_x^\ell \comm{\cH}{U}\partial_x^m V}_{L^p}\leq C\norm{ \partial_x^{\ell+m}U}_{L^\infty}\|V\|_{L^p}, && p\in(1,\infty), && \ell,m\in\mathbb{N}.
\end{align}
Using \eqref{commutatorH} and the Sobolev embedding
$$
\|\partial_x \Lambda^{-1}u \|_{L^\infty}\leq C\|\Lambda^{-1}u \|_{H^2}\leq C\|u\|_{H^1},
$$ 
we find that
$$
NL_2\leq C\|u\|_{H^3}\|u\|_{H^1}^2\leq \sigma\|u\|_{H^3}^2+C\|u\|_{H^1}^4.
$$
Similarly,
$$
NL_3\leq C\|u\|_{H^3}\|u\|_{H^2}\|u\|_{H^1}\leq \sigma\|u\|_{H^3}^2+C\|u\|_{H^2}^4,
$$
$$
NL_4\leq C\|u\|_{H^3}\|u\|_{H^2}\|u\|_{H^1}\leq \sigma\|u\|_{H^3}^2+C\|u\|_{H^2}^4.
$$
Integrating by parts in $NL_5$ and using the regularizing effect from $\mathcal{N}$, we can obtain that
$$
NL_5\leq C\|u\|_{H^3}\left(\|\partial_x u\|_{L^4}^2+\|u\|_{L^\infty}\|u\|_{H^2}\right).
$$
Using the Sobolev embeddings
$$
\|g\|_{L^4}\leq C\|g\|_{H^{0.25}},
$$
and
$$
\|g\|_{L^\infty}\leq C\|g\|_{H^{1}},
$$
we find that
$$
NL_5\leq \sigma\|u\|_{H^3}^2+C\|u\|_{H^2}^4.
$$
Integrating by parts and using the previous ideas we can estimate the last nonlinear contribution as
$$
NL_6\leq C\|u\|_{H^3}\|u\|_{H^2}^2\leq \sigma\|u\|_{H^3}^2+C\|u\|_{H^2}^4.
$$
Taking now $0<\sigma\ll$ small enough we can ensure that
\begin{align}\label{eq:energy}
\frac{\dd}{\dd t}\|u\|_{H^2}^2+\frac{\delta}{2}\|\Lambda^{3}u\|_{L^2}^2&\leq C\|u\|_{H^2}^2+C\|u\|_{H^2}^4,
\end{align}
which ensures the existence of a uniform time $T^*$ such that
$$
u\in C([0,T^*),H^2)\cap L^2(0,T^*;H^3).
$$
The local existence of solution will follow now from a standard application of Picard's theorem to a sequence of approximate problems. At this level of regularity, the uniqueness of such local strong solution can be easily obtained from a standard contradiction argument that we skip for the sake of brevity. \\

 The rest of this section is devoted to the global existence of solution for small initial data. In order to do that, we define the modified energy
$$
\nnorm{u}_{T} =  \ e^{\delta/2 T}\max_{t'\in\bra{0, T}} \set{ \norm{u\pare{t'}}_{A^0}} 
+\norm{u\pare{t'}}_{H^{2}}  .
$$
Then, our goal is to conclude the polynomial inequality
$$
\nnorm{u}_T\leq \cC_0\pare{f_0}+ P (\nnorm{u}_T).
$$

\subsection{The linear semigroup} We consider the linear nonhomogeneous problem
\begin{equation}\label{eq:DWW_asint2b}
2f_{t}-\mathcal{N}(f_{x}+2\delta \partial_x^2f+\cH f)+\beta \mathcal{P}\cH f_{xx}-\beta\delta \mathcal{P}\Lambda f_{xx}-\delta^2\mathcal{P}\partial_x^3 f+\delta^3\mathcal{P}\partial_x^4 f =F,
\end{equation}
where $F$ is the forcing. This linear equation can then be written as
\begin{align*}
2u_{t}+\mathscr{L} u =F,
&& u = \Lambda f, 
\end{align*}
with
$$
\widehat{\mathscr{L}u}=\lambda(k)\hat{u}(k)
$$
and
\begin{multline*}
\lambda(k)=-\frac{1-\delta ik}{1+\delta^2k^2}(ik-2\delta|k|^2 -i\text{sgn}(k)) \\
+\frac{\beta ik|k|}{1+\delta^2k^2}+\frac{\beta\delta|k|^3}{1+\delta^2k^2} +\frac{i\delta^2k^3}{1+\delta^2k^2}+\frac{\delta^3 |k|^4}{1+\delta^2k^2}
\end{multline*}
Then, we have that the homogeneous problem satisfies
\begin{align*}
\hat{u}(k,t)= & \  \hat{u}_0(k)e^{-\lambda(k)t}, \\
\av{\hat{u}_0(k)e^{-\lambda(k)t}} \leq & \  e^{-\delta t}\av{\hat{u}_0(k)}, 
\end{align*}
which in turn implies that
\begin{equation} \label{eq:control_uLb}
\ e^{\delta T}\max_{t'\in\bra{0, T}} \norm{u(t')}_{A^0} \leq \norm{u_0}_{A^{0}}. 
\end{equation}
Equivalently, we have that
\begin{equation}
\label{linear_operator_continuityb}
\|e^{-t \mathscr{L}}\|_{A^{0}\mapsto A^0}\leq  e^{-\delta t}.
\end{equation}

\subsection{Decay in the low regularity space}
Using Duhamel's principle, we can write the mild formulation of our problem as
$$
\hat{u}(k,t)=e^{-\lambda(k) t/2}\hat{u}_0(k)+e^{-\lambda(k) t/2}\int_0^t e^{\lambda(k) s/2}\hat{F}(k,s)ds
$$
with
\begin{multline*}
F=-\frac{\mathcal{N}}{2}\bigg{\lbrace}2uu_x+\Lambda\comm{\cH}{\Lambda^{-1}u }u +\beta\Lambda\comm{\cH}{\Lambda^{-1}u}\Lambda^2 u \\
-\delta\Lambda\comm{\cH}{u}u_x+\delta\partial_x\left(uu_x \right) +\delta\Lambda\comm{\partial_x^2}{\Lambda^{-1}u} u
\bigg{\rbrace}.
\end{multline*}
We observe that
\begin{align}
\widehat{\left(\comm{\mathcal{H}}{a} b\right)}\pare{n}&= \sum_n -i\left(\sgn k -\sgn\pare{k-n} \right) \ \hat{a}(n) \hat{b}(k-n), 
\end{align}
from where
\begin{align*}
0\leq |k|\leq |n|,
\end{align*}
so that this commutator does not vanish. A consequence of the above monotonicity relation is that
\begin{equation*}
\av{n-k}\leq \av{n}. 
\end{equation*}
 This implies that the above bilinear form presents a nontrivial commutation which allows to commute any derivative acting on the entire bilinear form as a differential operator acting onto $ a $ only. In a similar fashion, we find that
\begin{align}
\av{\widehat{\comm{\mathcal{H}}{\Lambda^{-1}u}\Lambda^2 u}\pare{k} }&= \av{\sum_n -i\frac{|k-n|^2}{\av{n}} \ \hat{u}(n) \hat{u}(k-n)\left(\sgn k -\sgn\pare{k-n} \right)}\nonumber\\
&\leqslant C \sum_n |n|\av{\hat{u}(n)}\ \av{\hat{u}(k-n)}.\label{commFourier}
\end{align}

Using that $A^0$ is an algebra, the fact that $\mathcal{N}$ gains one derivative, \eqref{commFourier} and Sobolev embedding, we find the estimate
$$
\|F\|_{A^0}\leq C \pare{ \|u\|_{A^0}^2+\|\Lambda^{-1}u\|_{A^0}\|u\|_{A^0}+\|u\|_{A^0}\|\Lambda u\|_{A^0} }.
$$
In particular  using
$$
\|u\|_{A^1}\leq C\|u\|_{A^{1.25}}^{4/5}\|u\|_{A^{0}}^{1/5}\leq C\|u\|_{H^{2}}^{4/5}\|u\|_{A^{0}}^{1/5},
$$
we find that
$$
\|F\|_{A^0}\leq C\|u\|_{A^0}^{6/5}\pare{ \|u\|_{A^0}^{4/5}+\|u\|_{H^{2}}^{4/5} }\leq Ce^{-3\delta/5 t}\nnorm{u}_{T}^2.
$$
As a consequence, we conclude that
\begin{multline}
\label{eq:loworder2}
 e^{(\delta/2) t}\max_{t'\in\bra{0, t}} \set{ \norm{u\pare{t'}}_{A^0}} \leq  C\|u\pare{t'}\|_{A^{0}}+C\nnorm{u\pare{t'}}_T^2\int_0^t e^{-(\delta/10) t'}\dd t'\\
 \leq  C\|u\pare{t'}\|_{A^{0}}+C\nnorm{u\pare{t'}}_T^2 , 
\end{multline}
this concludes the low-regularity estimates.

\subsection{Boundedness in the high regularity space}
To achieve the required estimate, we have to perform a finer analysis of the nonlinearity. In particular, we need to remove the term 
$$
\|u\|_{H^2}^2
$$
from the right hand side of \eqref{eq:energy}. In order to do that, we compute
\begin{align*}
NL_1&=-2\int_{\mathbb{S}^1}\mathcal{N}(uu_x)\partial_x^4 u \dd x,\\
&=2\int_{\mathbb{S}^1}\partial_x^2\mathcal{N}(uu_x)\partial_x^3 u \dd x,\\
&=2\int_{\mathbb{S}^1}\partial_x(1-\delta \partial_x)\mathcal{P}^{1/2}(uu_x)\mathcal{P}^{1/2}\partial_x^3 u \dd x,\\
&\leq C\|u\|_{H^1}\|u\|_{H^2}\|\mathcal{P}^{1/2}\partial_x^3 u \|_{L^2}.
\end{align*}
Similarly, invoking \eqref{commutatorH}, we find that
\begin{align*}
NL_2&=\int_{\mathbb{S}^1}\partial_x(1-\delta \partial_x)\mathcal{P}^{1/2}\Lambda\comm{\cH}{f}u \mathcal{P}^{1/2}\partial_x^3 u \dd x,\\
&\leq C\|u\|_{L^2}\|u\|_{H^{1.75}}\|\mathcal{P}^{1/2}\partial_x^3 u \|_{L^2},
\end{align*}
\begin{align*}
NL_3&=\beta\int_{\mathbb{S}^1}\partial_x(1-\delta \partial_x)\mathcal{P}^{1/2}\Lambda\comm{\cH}{f}\Lambda^2 u\mathcal{P}^{1/2}\partial_x^3 u \dd x,\\
&\leq C\|\Lambda^2 u\|_{L^2}\|u\|_{H^{1.75}}\|\mathcal{P}^{1/2}\partial_x^3 u \|_{L^2}.
\end{align*}
For the term $NL_4$ we compute as follows

\begin{align*}
NL_4&=\delta\int_{\mathbb{S}^1}(1-\delta \partial_x)\mathcal{P}^{1/2}\Lambda\comm{\cH}{u}\partial_x u \mathcal{P}^{1/2}\Lambda^4 u \dd x,\\
&\leq C\|\partial_x u\|_{L^2}\|u\|_{H^{1.75}}\|\mathcal{P}^{1/2}\Lambda^4 u \|_{L^2}.
\end{align*}

Similarly, since $ \cP^{1/2} \partial_x \pare{u u_x} = m_0\pare{D} \partial_x\pare{u^2} $ with $ m_0 $ a Fourier multiplier of order zero and using the classical fact that the space $ H^s\cap L^\infty, \ s\in\bR $ is a Banach algebra we can argue that

\begin{align*}
NL_5&=-\delta\int_{\mathbb{S}^1}(1-\delta \partial_x)\mathcal{P}^{1/2}\partial_x\left(uu_x \right)\mathcal{P}^{1/2}\Lambda^4 u \dd x,\\
&\leq C(\|u\|_{H^2}\|u\|_{A^{0}}+\|\partial_x u\|_{L^4}^2)\|\mathcal{P}^{1/2}\Lambda^4 u \|_{L^2},\\
&\leq C\|u\|_{H^2}\|u\|_{A^{0}}\|\mathcal{P}^{1/2}\Lambda^4 u \|_{L^2},
\end{align*}
and
\begin{align*}
NL_6&=-\delta\int_{\mathbb{S}^1}(1-\delta \partial_x)\mathcal{P}^{1/2}\Lambda\comm{\partial_x^2}{f} u \mathcal{P}^{1/2}\Lambda^4 u \dd x,\\
&\leq C\|u\|_{H^2}\|u\|_{A^{0}}\|\mathcal{P}^{1/2}\Lambda^4 u \|_{L^2},
\end{align*}
where we have used the inequality
$$
\|\partial_x u\|_{L^4}^2\leq C\|u\|_{H^2}\|u\|_{A^{0}}.
$$
As a consequence, using Young's inequality, we can find the inequality
\begin{align}\label{eq:energyv2}
\frac{\dd}{\dd t}\|u\|_{H^2}^2+\frac{\delta}{2}\|\mathcal{P}^{1/2}\Lambda^4 u \|_{L^2}^2&\leq C\nnorm{u\pare{t'}}_T^3\|u\|_{H^2}e^{-\delta/16 t}.
\end{align}

\subsection{Finishing the proof of Theorem \ref{theorem2}} 
Collecting \eqref{eq:loworder2}, \eqref{eq:energyv2} we conclude the polynomial bound
$$
\nnorm{u}_T\leq \cC_0\pare{u_0}+ \nnorm{u}_T^2.
$$
From here we can finish the argument as in the proof of Theorem \ref{theorem} and we obtain that that the solution is global.

\section*{Acknowledgments}
The research of S.S. is supported by the European Research Council through the Starting Grant project H2020-EU.1.1.-639227. R.G-B was supported by the project ”Mathematical Analysis of Fluids and Applications” with reference PID2019-109348GA-I00/AEI/ 10.13039/501100011033 and acronym ``MAFyA” funded by Agencia Estatal de Investigaci\'on and the Ministerio de Ciencia, Innovacion y Universidades (MICIU). Project supported by a 2021 Leonardo Grant for Researchers and Cultural Creators, BBVA Foundation. The BBVA Foundation accepts no responsability for the opinions, statements and contents included in the project and/or the results thereof, which are entirely the responsability of the authors.
 
\providecommand{\bysame}{\leavevmode\hbox to3em{\hrulefill}\thinspace}
\providecommand{\MR}{\relax\ifhmode\unskip\space\fi MR }
\providecommand{\MRhref}[2]{%
  \href{http://www.ams.org/mathscinet-getitem?mr=#1}{#2}
}
\providecommand{\href}[2]{#2}

\end{document}